# "Audacity or Precision":
# The Paradoxes of Henri Villat's Fluid Mechanics in Interwar France

David Aubin
May 2010

**Abstract:**
In Interwar France, Henri Villat became the true leader of theoretical researches on fluid mechanics. Most of his original work was done before the First Word War; it was highly theoretical and its applicability was questioned. After having organized the first post–WWI International Congress of Mathematicians in 1920, Villat became the editor of the famous *Journal de mathématiques pure et appliqués* and the director of the influential book series "Mémorial des sciences mathématiques." From 1929 on, he held the fluid mechanics chair established by the Air Ministry at the Sorbonne in Paris and was heading the government's critical effort invested in fluid mechanics. However, while both his wake theory and his turbulence theory seemingly had little success outside France or in the aeronautical industry (except in the eyes of his students), applied mathematics was despised by the loud generation of Bourbaki mathematicians coming of age in the mid 1930s. How are we to understand the contrasted assessments one can make of Villat's place in the history of fluid mechanics?

## Introduction

The title of this paper may call for an explanation. On 16 June, 1911, a young mathematician named Georges Léry complimented his friend Henri Villat on the doctoral thesis he had justr completed: "One does not know whether to admire the most in your work the audacity or the precision of your reasoning."[1] At the time, Léry was a mathematics teacher at the *lycée* of Reims. He had befriended Villat—who would in the interwar period become the leader of fluid mechanics research in France—at the Rollin College from where they proceeded to the prestigious École normale supérieure (ÉNS), Paris, in 1899.[2] Now, to me, "audacity or precision" stands for one of the fundamental conundrums faced by theoreticians interested by the study of fluid flows. While there often was only one type of precision most specialists in fluid mechanics would have dreamed to achieve, namely deriving phenomenologically plausible conclusions from fundamental equations, there were several options—when, as it is so often the case in this domain, this dream proved to be unattainable—to display one's theoretical audacity. Should one uncompromisingly stick to hypotheses closely mirroring reality—study compressible viscous fluids, for example? Or should one instead focus on experimental results and explain them by only loosely appealing to fundamental fluid equations? Or else, might one be entitled to make daring hypotheses about the fluid behaviour, for example introducing discontinuities in the flow, in the hope of closing the gap between mathematical deductions and experimental results? In other words, what was fluid mechanics: physics, engineering, or mathematics?[3]

---

[1] "On ne sait ce qu'il faut le plus admirer dans ton travail, l'audace ou la précision des raisonnements. » Léry to Villat (16 june, 1911). The 8 boxes of the "Fonds Villat" are available at the Archives of the Academy of Sciences, call number 61J. Most of the correspondence I have used is to be found in box 5. I will hereafter refer to it as Villat Correspondence.

[2] On Léry's relationship with Villat, see (Hadamard 1916). Other biographical information on Villat's life and work can be found in (Leray 1973, Roy 1972, Jacob 1979). Founded as a branch of the famous Sainte-Barbe College, the Rollin College is now called the lycée Jacques Decours and is located near Montmartre in Paris.

In the case of Villat's work, both audacity and precision were seemingly achieved. As for paradoxes, Villat's career can easily placed under that banner. To start with, most of his work was concerned with the so-called D'Alembert paradox, the fact that the steady motion of a perfect fluid past an obstacle produces no pressure on the obstacle, so that birds—and airplanes—cannot fly! Using the methods developed by Hermann von Helmholtz and Gustav Kirchhoff some fifty year before, and just perfected by Tullio Levi-Civita, Villat produced something that struck contemporaries as a mathematical monument of classical beauty. In his analysis, Villat introduced discontinuous surfaces that needed to be extended indefinitely, which was named the Brillouin paradox.

But, to me, the greatest paradox was the professional success of Villat himself. It was said that Villat's main contribution lay in making Levi-Civita's work "applicable" (Jacob 1979). But, as Maurice Roy diplomatically explained, if his mentor's work "remain[ed] applicable" in theory, it "ultimately ran into the fact that the field of practical use of wake theory was constantly reduced by the search for obstacles and flow regimes reducing wakes—a search soon to become dominant in applications of fluid mechanics" (Roy 1972). But even while his main contributions to fluid mechanics appeared less and less relevent to practical applications, Villat reached a central position on the French mathematical scene. Here was a man who, at the age of 40, was put in charge of the material organization of the Strasbourg International Congress of Mathematicians in 1920, then placed in a key editorial position at both the *Journal de mathématiques pures et appliquées* and the *Bulletin des sciences mathématiques*, who was given a prestigious and well-endowed chair of fluid mechanics at the Sorbonne in Paris in 1929 and who directed a large number of doctoral theses. In short, building on his unapplicable theory of mixed mathematics at a time when applications seemed crucial in public discourses,Villat was able to become one of the "bosses" of French mathematics (Gispert et Leloup 2009, Leloup 2009).

Although one may say that Villat was tremendously important in the institutional landscape of interwar French mathematics, his theoretical and mathematical work in fluid mechanics seems to have had limited impact at the time and since. As the historian Hashimoto has written about wake theory (which he called stream theory): "Such a deep commitment to the discontinuous stream theory surprises modern readers. With so uncertain a physical basis, the construction of such a sophisticated mathematical theory may even appear absurd" (Hashimoto 1990, 149). Villat's hope was that wake theory provided the best approximation of the zero-viscosity limit for fluid flows. Although this has sometimes been forgotten, he was not alone in thinking this. In the British Advisory Committee for Aeronautics, George Greenhill was a strong adherent to that theory. Stefan Drzeviecki, who translated N. Zhukovsky's *Aérodynamique* (Paris, 1916), wrote that while Levi-Cività introduced an arbitrary function helping to determine the current flowing past an obstacle: "Villat, whose work is more complete, shows the manner to determine this arbitrary function on the condition that one knows the form of the contour along which the flow occurs as well as the points on this contour where the string of fluid is detached from it."[4] Reviewing the state of wake theory in 1930, the engineer Maurice Biot wrote that it had achieved its definitive state of development through Villat's contributions. "Unfortunately," he added, "theoretical results only very roughly account for data. Experiment has shown that in many cases […] they do not conform to the regime that is established" (Biot 1930, 237). Disagreement between theory and experiment however was no impediment to Villat's professional rise. Why was it so?

---

[3] For a general history of fluid mechanics, see (Darrigol 2005). I have discussed Marcel Brillouin's, Henri Poincaré's, and Henri Bénard's respective views on fluid mechanics around 1900 in (Aubin, 'The Memory of Life Itself': Bénard's Cells and the Cinematograph of Self-Organization 2008).
[4] Quoted and trans. in (Hashimoto 1990, 150).

In the following, to try to answer this question, I take up a point recently made about the effect of World War I on French mathematics (Goldstein 2009). Among French mathematicians, as for many segments of French society, the Great War represented a trauma whose full contours we are barely starting to discern properly.[5] As Catherine Goldstein has argued, we now understand that several conceptions of mathematicians' engagement in service of the collective then coexisted and, in particular, that the intense production of original research could sometimes be side-stepped in favor of other forms of commitment. To understand Villat's trajectory, as we shall see, his obvious organizing and pedagogical skills were compounded with a specific ideology of applied mathematics. His approach to the turbulence problem entered in resonance with discourses emanating from wartime experiences.[6]

To make my claim clear, I review Villat's work before World War I, focusing on his doctoral thesis and the controversies it brought about with contemporary French mathematicians such as Marcel Brillouin, Paul Painlevé and Pierre Duhem. I then turn to Villat's activities in Strasbourg after the war and in Paris in the late 1920s. Finally, Icome back to wartime and postwar discourses about mathematics and science to show that the brand of fluid mechanics brandished by Villat might be exactly what was needed at the time: a classic branch of French analysis which although developed out of disinterested motives was not completely detached from reality. In other words, Villat's mathematical theory of fluid flows could be criticized neither for being purely abstract nor for being merely practical. In the context of postwar France, this, for many, was precisely the kind of mathematics one needed to develop.

## Fluid Resistance and the Airplane Problem

After his studies at the ÉNS, Villat assumed the position of mathematics teacher in a secondary school, the *lycée* of Caen, Normandy, in October 1902. His mathematical talent and higher ambition were apparently quite obvious to his entourage and already in 1906 he was asked to lecture at the University of Caen. The seed of his interest in fluid mechanics is unknown, but as the old Sorbonne hydrodynamics specialist Joseph Boussinesq wrote to him, this might have to do with the airplane problem.[7]

Aviation obviously was a hot topic in the first decade of the 20th century. But, it was unclear who would be in the best position to come up with a satisfactory theory of flight. Around 1910, the intellectual and social origins of the pioneers of aviation were quite diverse. The first convincing flight that had taken place in France, on 13 January, 1908, at Issy-les-Moulineaux, was the combined effort or three graduates of the École des Beaux-Arts in Lyons. The airplane was built by Gabriel Voisin (who had an architect's training), equipped with an engine built by Léon Levavasseur (trained as a sculptor and an autodidact in mechanics) and piloted by Henry Farman (who was trained as a painter and was a cyclist turned automobilist). Among French pioneers people with an engineer's training, like Louis Blériot, Louis Breguet or Ferdinand Ferber, in fact predominated. An artillery officer trained at the École polytechnique, Captain Ferber clearly stood out due to the scientific considerations and the careful experimental method he applied to the problem of flying as early as 1898. An historian of aviation assessed that the influence of this mixed bunch on the French aeronautical industry in "the first half of the 20th century was uninterrupted and considerable," imprinting an entrepreneurial "style" on this industry, which as result never became "an activity directly linked to science" (Chadeau 1985, 283).

---

[5] A research group at the Institut de mathématique de Jussieu has taken up the study of mathematics and mathematicians around WWI. Some of its findings have been published in a special issue of the *Revue d'histoire des sciences* (Gispert et Leloup 2009, Goldstein 2009, Beaulieu 2009) and in an edition of Vito Volterra's correspondence with French mathematicians (Mazliak et Tazzioli 2010), and will be in two edited books in preparation (**refs\*\*\***). See also D. Aubin, "L'élite sous la mitraille: mathématiciens, mémoire et oubli de la Grande Guerre," in preparation. These work will bring useful counterpoints to earlier assessments of the effect of WWI on French scientists (Roussel 1989, Prochasson et Rassmussen 1996, Aubin et Bret 2003).

[6] This point obviously owes to Forman's historiographical point about the influence of Weimar culture on the emergence of quantum mechanics (Forman 1971). This thesis has, as it is well known, been much criticized, a point I will address in the conclusion of this article.

[7] For a mention of airplanes, see Boussinesq to Villat (11 May, 1910), Villat Correspondence.

Neither was academic science completely out of the picture. Since 1902, the Academy of Sciences had put together an Aeronautics Commission that included a dozen of scientists from every specialty: astronomers, physicists, chemists, engineers, mechanicians, mathematicians, and even physiologists including the pioneer of cinematography Étienne-Jules Marey who had a long-standing interest in the flight of birds.[8] From this variety of approaches, many concurrent scientific takes on mechanical flight were developed. But, of all these approaches, the mathematical approach insisting on using fundamental equations of fluid motion to derive the principles of flight was far from being the most popular. It became a running joke that the theoreticians had rather concluded that airplanes could not fly. In 1904–1905, the ballistician Emmanuel Vallier had indeed published a series of careful theoretical studies in which he rejected the practical possibility of flying with an airplane because he had estimated that the horizontal velocity needed for producing an adequate lift was too high.[9] Barely two years later, men with a practical bent would triumphantly claimed that this mathematical problem had been solved by the airplane itself: "luckily, M. Santos-Dumont's trials at Bagatelle has refuted Colonel Vallier's pessimistic conclusion" (Armengaud 1907, 867).[10]

<Insert Figure 1 here>

When Orville and Wilbur Wright exhibited over France their *Flyer*, the mathematician and politician Painlevé was one of the first to embark, on October 10, 1908. Flying for over an hour above Auvours, he and Wilbur in fact held the record for the longer passenger flight for some time afterwards (Figure 1).[11] As a member of Parliament, Painlevé quickly realized the military and commercial implications of the new invention. As early as February 1909, he submitted to the Parliament a project for a national aviation laboratory. Although this never saw the light of day, an Aerotechnical Institute [*Institut d'aérotechnique*] was soon set up in Saint-Cyr near Paris and attached to the university of Paris (Figure 2). While the physicist Charles Maurain was appointed as director of Institute, a new professorial chair in aviation was offered to another physicist Lucien Marchis (Maurain 1913, Maurain 1914). A special technical school for aeronautics [*École supérieure d'aéronautique et des constructions mécanique*] was inaugurated in November 1909. Much of these institutions were privately funded by industrialists such as Basil Zaharof and Henry Deutsch de la Meurthe. To present synthetic views of the scientific and technological development of aeronautics, a new journal was launched in 1910, *La Technique aéronautique,* which counted scientists, inventors, industrialists and military officers among its editorial board (Anizan 2006, 132-134).

---

[8] The archives of the Commission are kept at the Academy of Sciences in Paris (DG 47) and to my knowledge have never been systematically studied. On Marey's fluid mechanical experiments, see (Didi-Hubermann et Mannoni 2004).
[9] Dynamique de l'aéroplane, *Revue de mécanique* (1904), p. 5, 101, 342, 539 and (1905), 125, 237, 323. ***
[10] From September to November 1906, Alberto Santos-Dumont flew an airplane with a Levavasseur engine for as much as 220 meters. For an earlier skeptical review of Vallier's work, see (Soreau 1905). See also the anonymous report of the Société française de physique meeting on 8 March, 1908, in the *Revue générale des sciences pures et appliqués* 19, (1910),416 as well as (Lecornu 1909).
[11] See (Anonymous 1909). On Painlevé and aviation, see (Fontanon et Franck 2006, 41-56, Fontanon 1998, Anizan 2006, 121ff). See also Paul Painlevé and Émile Borel, *L'aviation* (Paris: F. Alcan, 1910).

< Insert Figure 2 here. >

As a mathematician, Painlevé immediately transformed the flying machine into a special problem for fluid mechanicians: "The motion of a solid in a fluid, air or water for example, can be studied mathematically when the viscosity of the fluid, the roughness of the edges and shocks are neglected" (Painlevé, L'aéroplane 1910, i). But this was a complex problem and Painlevé's scientific work as well as his political action were geared, not toward a mathematical solution, but rather toward a better integration of theoretical, experimental, and technological research. The problem, he went on:

> leads to extremely complex computations and, this notwithstanding, neglected factors play a considerable role in phenomena that are of interest for aviation. So, if one wants to take these factors into account, the problem becomes intractable. One has to have recourse to experimentation (ibid.).

It is therefore in a context where excitement was great for the scientific study of aerial travel but confidence low in mathematical methods to solve it that Villat breached the main problem of fluid mechanics. Although fruitful discussions between theoreticians and inventors were hard to established, they could meet around the problem of the resistance exerted by air on a solid body moving through it. Resistance indeed was a quantity that seemed both computable and observable through experimentation. Moreover, the observed value of air resistance in the airplane was estimated to be roughly ten times lower than expected (Lecornu 1909, 36). The classic expression for air resistance simply assumed that resistance was proportional to the square of the velocity (at least for small velocities) and to the cross section of the obstacle: $F = KSV^2$. The coefficient $K$ depended on the shape of the body and teams headed by Marey and by Gustave Eiffel had performed many experiments to determine its value for different objects (Renard 1909, 83-84).[12] Villat's initial success as a theoretician lies in the fact that he was the first to be able to derive a value for air resistance from (Euler) equations of motion for a very large class of obstacles of different shapes.

## Boussinesq's and Brillouin's Mentorship

Why, contrary to Painlevé and most people involved in aeronautics at the time, would a young mathematics professor in Caen think that air resistance could be computed from first principles is not clear. But when Villat set out to work on the problem of the resistance exerted by a solid in a fluid flow, it did it with such ardor that he apparently wrote the draft of six articles in the course of less than a year. In a manner that seems premonitory in view of his later role in publishing, Léry enthusiastically wrote: "This is phenomenal. Soon there will not be enough periodicals to welcome them all and you will need to establish a new one."[13] But it is not easy to determine precisely what Villat's approach was at the time and from where he drew his inspiration. The first trace we have of his tackling the problem comes from the he received from Boussinesq, dated 11 May, 1910, mentioned earlier.

---

[12] On Marey, see references given in note 8. On Eiffel, see (Lee, Fontanon et Fruman 2005) and Claudine Fontanon's article in this collection and (Eiffel 1907).

Despite his old age Boussinesq still was the active leader of French theoretical fluid mechanics, but seemed not concerned too much with aviation. He indeed was less than enthusiastic about Villat's request to see him to discuss the problem of air resistance applied to airplane, emphasizing that a trip to Paris for that purpose was quite unwarranted. "What little things I would tell you, I'll put in this letter."[14] In this letter, Boussinesq noted that before dealing with the compressible case corresponding to high velocities, one should try to solve the problem for perfect fluid. Villat would stick to this first approximation for most of his later career as a fluid mechanician. According to Boussinesq, the best source on this problem was the massive memoir he had himself published in 1887 on the resistance of fluids.[15] Still, he warned, Villat would find extremely little in there. Boussinesq seemed as doubtful as any about the success of a theoretical approach to the problem, insisting on the need "almost blindly to do an infinity of observations, always to be repeated whenever the shape and even perhaps the dimension of the engine would change."[16] The only case were "*haute analyse*" provided some answer, he added—the slow motion where nonlinear terms in the equation were eliminated—was of no interest to airmen.[17]

Boussinesq nonetheless provided Villat with a precious reference to Italian literature which he had not read due to the language in which they were written. In 1909–1910, several papers published by Emilio Almansi and a work by Levi-Civita cited herein had nonetheless seemed noteworthy to the Sorbonne professor. Considering the lack of reference to the Italian literature in Villat's first published article on fluid mechanics dated June 1910(Villat 1910a), we may infer that this directed him to Levi-Civita's work which played a crucial part in Villat's own approach to the problem of air resistance. Boussinesq however doubted that this approach would lead to great advances since the Italian authors had restricted their attention to perfect fluids when it was obvious that even at small speed, vortices in the wake on the immersed body played a fundamental role.

> But you are no doubt a young man, full of energy and desirous to use it. You would not be able to do so in a more useful manner than on such a momentous problem to which a solution has become so needed to fulfill the wishes shown by your restless generation.[18]

---

[13] « c'est phénoménal. Il n'y aura bientôt plus assez de périodiques pour les accueillir et tu seras obligés d'en fonder un ». Léry to Villat (15 March, 1910). Villat Correspondence.

[14] « Le peu que je vous dirais, je le mettrai dans cette lettre. » Boussines to Villat (11 May, 1910), Villat Correspondence. On Boussinesq, see (Darrigol 2005)

[15] Boussinesq, "Sur la résistance des fluides," *Mémoires de l'Académie des Sciences*, 14 (1887).

[16] « <u>ce tout</u> est <u>extrêmement peu de chose, et ne dispensera nullement</u> de faire, à peu près en aveugle, une infinité d'observations, toujours à recommencer pour peu que changent les formes, et, peut-être même, les dimensions des appareils. » Boussinesq to Villat (11 May, 1910), Villat Correspondence. Emphasis in the original.

[17] As Boussinesq pointed out this case had been studied in his own *Cours de physique mathématique*, vol. 2 (Paris: Gauthier-Villars, 1903), 199-264.

[18] « Mais vous êtes, sans doute encore jeune, plein d'ardeur et désireux de la dépenser. Vous ne sauriez le faire plus utilement que sur un problème si actuel et, on peut le dire, devenu si nécessaire à résoudre pour répondre aux besoins que manifeste notre remuante génération. » Boussinesq to Villat (11 May, 1910), Villat Correspondence.

Villat immediately sent his reply with a more detailed explanation that has not been found, but Boussinesq's reaction to it is instructive, even if he admitted that he had not studied Villat's work as carefully as he should have. In the permanent case, the old professor considered how the lines separating the region where the fluid remained perfect from the turbulent one closed back in the wake of the body and mentioned the relatively simple solution found by Gustav Kirchhoff where those lines remained parallel at infinity, emphasizing: "at bottom, [this solution has] almost nothing to do with our problem."[19] Again, this was pointing in the direction of the approach adopted by Villat in his thesis.

Villat had poked Boussinesq's interest. He now acknowledged that Villat's work might be sufficient for a doctoral thesis, especially in this field where expectations should be kept low:

> On such a difficult question, in case I would have to examine your work, I would not feel the right to be demanding regarding originality; and a few personal ideas along the lines of those you have (much too briefly) submitted to me could *à la rigueur* be sufficient, even if they seem doubtful, provided they would be supplemented with a detailed historical account of more or less recent work like for example, in Italy, those of Almansi, Lev-Civita and also Louis de Marchi (if I recall correctly) [sic]. Such an exposition would be helpful and would provide you with the opportunity of showing the most serious qualities of the scientist and researcher. So, don't be discouraged.[20]

At this juncture, Villat might have realized that not much scientific exchange was to be had with Boussinesq. Over the summer 1910, he started another scientific correspondence, this time with Marcel Brillouin. As general and mathematical physics professor at the Collège de France he had lectured on fluid mechanics some years earlier (Aubin 2008), and again in 1908–1909 in the wake of the first flights, on the resistance of fluids and gases to the motion of solid bodies through it. This letter exchange would prove much more satisfying to the aspiring fluid mechanician.

<Insert Figure 3 here>

On July 27,1910, Brillouin thus wrote to Villat that his thesis project seemed very interesting to him, even if he were not to add anything new to Levi-Civita's work. Writing that he had already achieved incomplete results on this topic, Brillouin incisively insisted on the aspects of the Italian mathematician's work that required some extensions. Villat should focus on the possible determination of the function $\Omega$ in terms of the geometrically defined shape of the obstacle—"*Voilà le but.*" Any solution of the problem for a specific surface would be a definite progress:

> if one could express $\zeta$ in terms of $\Omega$ (instead of $\Omega$ in terms of $\zeta$), the computation of $z$ would be, except for a finite or infinite number of computable terms) to integrals of the form

---

[19] « au fond presque étrangère à notre problème. » Boussinesq to Villat (1 June, 1910), Villat Correspondence.
[20] « Je ne m'étais pas, d'abord, douté que vous voulussiez faire une thèse de doctorat. Dans une question si difficile, et supposé que je fusse examinateur, je ne me reconnaîtrais pas le droit d'être très exigeant sous le rapport de l'originalité ; et quelques idées personnelles, dans le genre de celle que vous me soumettez (beaucoup trop succinctement) pourraient à la rigueur suffire, parussent-elles même un peu douteuses, <u>pourvu qu'elles fussent accompagnées d'un historique assez nourri, bien clair, des travaux plus ou moins récents</u>, comme par exemple, en Italie, ceux d'Almansi, Levi-Civita et aussi Louis de Marchi (si je me souviens bien). Un pareil exposé rendrait des services et vous permettrait de montrer les qualités les plus sérieuses du savant et du chercheur. Ne vous découragez donc pas. » Boussinesq to Villat (1 June, 1910), Villat Correspondence.

$$\int e^{-qx}(1-\cos x)\frac{dx}{x}, \text{ and } \int e^{-qx}\sin x\frac{dx}{x}$$

whose tables do not yet exist, but that I have started to compute in view of this application.[21]

To discuss fully the difficult resolution of $\zeta = f(\Omega)$, Brillouin added, one would need to rely on recent knowledge about meromorphic functions by Jacques Hadamard, Émile Borel, and Gosta Mittag-Lefler, which he acknowledged to be lacking. Himself had worked out some special cases by numerical and graphical means (see Figure 3). "But, it seems to me that *you* would be able to undertake a classification of surfaces according to the form of the function $\zeta = f(\Omega)$, and from this angle tackle the problem of correlation between Ω and the shape of the surface."[22]

## Fluid Resistance and Wake Theory

To understand what is at stake in Brillouin's letter, it is now necessary to delve deeper in the mathematical content of Villat's doctoral thesis (Villat 1911). There, he positioned his work in continuity with Helmholtz's classic discontinuous wake theory first suggested as early as 1868. In this theory, the resistance of a perfect fluid to the motion of a solid body was modeled by a surface of discontinuity behind the body. This particular solution to D'Alembert's paradox had been explicitly solved in the case of rectilinear obstacles. As mentioned earlier Levi-Civita made a significant contribution in 1907 to the problem by using conformal transformations, a tool developed in the field of complex analysis.[23]

< Insert Figure 4 Figure 4 : **Representation of the obstacle in the *z*-plane corresponding to the physical motion of the fluid around the obstacle .**here. >

To assess Villat's contribution to Levi-Civita's work, let us examine the mathematical formalism they mobilized. Consider the motion of an incompressible fluid in the *xy*-plane in which a solid obstacle was immerged. Most supposed that the body was at rest and that the asymptotic horizontal component of the fluid velocity was constant and positive while its asymptotic vertical component was nil (Figure 4). The flow regime was supposed to be permanent and irrotational Let *O* be the point on the body where fluid lines split into two parts (this could be the pointy tip of a ship or a shell). This point was taken to be the origin of the (*x*,*y*) coordinate system. The region *(A)* of the plane where the fluid was moving with respect to the body was delimited by two lines composed of (1) the profile of the solid $\varpi_1 + \varpi_2$ and (2) two lines of discontinuity $\lambda_1$ and $\lambda_2$ where the velocity field suddenly jumped from one value to another.

To analyze this problem, Levi-Civita and Villat followed an idea well known since Lagrange and Cauchy consisting in expressing it as a problem in terms of functions of complex variables.[24] At a point (*x*,*y*) the velocity was defined by its two components *u* and *v* satisfying the incompressibility condition:

$$\frac{\partial u}{\partial x} + \frac{\partial v}{\partial y} = 0.$$

---

[21] « si on pouvait exprimer ζ en Ω (au lieu de Ω en ζ), on ramènerait le calcul de z, outre un nombre fini ou infini de termes calculables) à des intégrales de la forme [voir formules ci-dessus] dont les tables n'existent pas encore, mais que j'ai entrepris de faire calculer, en vue de cette application. » Brillouin to Villat (27 July, 1910). Villat Correspondance. These results which Brillouin had introduced in his lessons at the Collège de France in March and April 1909. The computations mentioned by Brillouin were in part done by Kannappel.
[22] « Aussi ne puis-je que soupçonner le degré de difficulté. Toutefois, il me semble que vous puissiez aborder un classement des formes de surface, d'après la forme de la fonction ζ = f(Ω), et attaquer de ce côté le problème de la corrélation entre Ω et la forme de la surface. » Brillouin to Villat (27 July, 1910). Villat Correspondance.
[23] For studies of fluid mechanics in Italy, see (Ricca 1996, Nastasi et Tazzioli 2005, 203-210).

The simplest way to express this condition was to suppose that there existed a function $\psi$ such that :

$$u = \frac{\partial \psi}{\partial y}, \qquad v = \frac{\partial \psi}{\partial x}.$$

Since irrotational motion was assumed here, the velocity field could moreover be expressed as the gradient of a function $\varphi$, called the velocity-potential:

$$u = \frac{\partial \varphi}{\partial x}, \qquad v = -\frac{\partial \varphi}{\partial y}.$$

Both functions $\psi$ and $\varphi$ being harmonic (that is, $\Delta\varphi = \Delta\psi = 0$), the problem was very easily transposed in the complex plane, and these relations could be reduced to a single one. Defining:

$$\begin{cases} z = x + iy \\ f = \varphi + i\psi \end{cases},$$

then *w* is a finite continuous complex function of *z* over the whole portion of the plane *z* occupied by the liquid in motion and:

$$w = u - iv = e^{-i\Omega} = \frac{df}{dz}.$$

< Insert Figure 5 here. >

This function $\Omega$ was the one to which Brillouin suggested Villat to pay close attention. The advantage of this approach was that, following an idea first developed by Kirchhoff (in *Vorlesung der Mechanik*, translated in French by Célestin Soutreaux in his thesis called *Sur une question d'hydrodynamique* [Paris, 1893]), it allowed mathematicians to apply conformal mapping theory. An analysis of the mathematical relation between *f* and *z* indicated that for any obstacle lying in the *z*-plane (Figure 4), there corresponded a simple situation in the *f*-plane (Figure 5). The region *(A)* in the *z*-plane was mapped by the conformal mapping $f(z)$ to the region *(B)* in the *f*-plane and the obstacle to a straight line which since the functions $\psi$ and $\varphi$ could be chosen up to a constant factor was chosen to be the semi-axis $\varphi = 0, \ \psi > 0$.

< Insert Figure 6 here. >

Levi-Civita further suggested that a proper change of variable allowed to map conformally the region *(B)* to a semi-annulus in a new complex plane $\zeta = \xi + i\eta, \ \eta > 0$ (Figure 6). But it was

---

[24] The following argument was standard at the beginning of the 20th century. It is for example explained in (Brillouin 1911, Levi-Cività 1907).

Villat who was first able to compute the equation corresponding to this change of variables in the most general case:

$$f = A\left[\frac{b}{3} + \wp\left(\frac{\omega}{i\pi}\log\zeta - \omega'\right)\right] + A(b-a)\log\left[-\frac{2b}{3} - \wp\left(\frac{\omega}{i\pi}\log\zeta - \omega'\right)\right] + C,$$

where $A$, $C$, $a$ and $b$ were constants that could be computed and $\wp(z)$ was known as the Weierstrass elliptic function with periods $\omega$ and $\omega'$. In this representation, the limits of the wake corresponded to the parts of the semi-annulus lying on the axis $O\xi$.

In 1907, Levi-Civita showed that conformal theory allowed the computation of the law of resistance corresponding to a body plunged in the fluid flow.

$$\mathbf{R} = R_x + iR_y = \frac{1}{2i}\oint_{|\zeta|=1} e^{i\Omega} df.$$

Exact determination was however not possible since his derivation depended on the unsolved integral. In his 1907 article, Levi-Civita used a power series to solve some simple cases including the so-called Babyleff problem where the obstacle was two straight lines at an angle.

Being explicit, Villat's solution was more "applicable," to use Jacob's term quoted above. Using the expression for $f$ derived above, Villat was able to find an explicit expression for $df$ as a very complicated function of $\zeta$. All these results then allowed him to express the resistance of the obstacle to the fluid flow as follows:

$$P = -\frac{A\omega}{2i\pi}(e_3 - e_1)(e_3 - e_1)^2$$
$$\times \int_0^{2\pi} \frac{e^{-T}(\cos\Theta + i\sin\Theta)\left[\wp\left(\frac{\omega}{\pi}\sigma_0\right) - \wp\left(\frac{\omega}{\pi}\sigma\right)\right]\wp'\left(\frac{\omega}{\pi}\sigma\right)}{\left[\wp\left(\frac{\omega}{\pi}\sigma_0\right) - e_3\right]\left[\wp\left(\frac{\omega}{\pi}\sigma\right) - e_2\right]\left[\wp\left(\frac{\omega}{\pi}\sigma\right) - e_3\right]^2} d\sigma.$$

Obviously, such a terrible-looking equation is impressive and quite a bit of further explanations would be needed to make its terms properly understandable. Let us rather focus here on the significance it may have had at the time.[25] Obviously, making full use of the modern theory of functions, Villat's approach had a mathematical sophistication which strongly contrasted with the elementary mathematical methods mobilized by "aerotechnicians" (Maurain 1914). But Villat's objective remained congruent with the latter's main concerns. As explained above, the problem of air resistance was important for aviation and the hope was that the air resistance coefficient $K$ in the empirical resistance law $F = KSV^2$ could be computed for any possible shape. Developing upon the above formula in his thesis, Villat was indeed able to compute explicit values of the resistance exerted by a few simple symmetrical shapes (Villat 1911, figs. 23, 24 and 25).

---

[25] The most detailed explanation can be found in Villat's thesis (Villat 1911, 228-230), as well as some of his later publications, for example (Villat 1920).

## Disappointment and Success

According to one of his students, Villat's last surviving letter from Boussinesq has had fateful consequences. In a necrology (Leray 1973), it is indeed reported that Villat had proven a version of the famous Kutta-Zhukovsky theorem according to which the lift of an infinitely long wing is proportional to the velocity flow around the wing (or in other words, the period of the velocity potential). But Boussinesq deemed faulty the result which were to become key to wing theory. Villat wished to see his result published in the widely distributed *CRAS* which could be quickly done if a member of the Academy (which Boussinesq was) recommended its insertion in the next issue. But, Boussinesq's negative answer by recalling a version of the D'Alembert paradox clearly shows that he had not followed Villat's reasoning:

> Your Note has reached me in time for its presentation tomorrow, Monday, in the last meeting to which I will be able to assist [before his going away for the holidays]. But I do not feel that I should present it; because the *only* case dealt with in its analytic part leads to the conclusion that $P_x = 0$, which is the precise opposite to your main and final conclusion concerning, so you say, the permanent *rotational* motion completely foreign to your formulas based on the use of a potential $\varphi$ implying an absence of vortices. A text that would properly warrant your conclusion should not at all contain the potential $\varphi$.[26]

Boussinesq moreover pointed out that for a permanent flow regime, it would be impossible for the obstacle to exert an impulsion on the fluid without giving rise to a corresponding momentum (while it is assumed to stay constant within the fluid) unless this momentum were dissipated at an infinite distance for the body, a hypothesis Villat discarded by neglecting surface integrals on a large circumference far from the body. Furthemore, Boussinesq went on, in case the fluid motion was not permanent, whether rotational or not, a wind gush or a sudden aspiration would surely lead to a non-vanishing horizontal component of the force $P_x$. "Nobody will learn anything for your Note."[27]

It is poignant to read the part of the letter where Boussinesq addressed questions of priority: "work published abroad during the time you are writing your thesis would not *to my eyes* strip you of this priority should I be asked to assess this thesis. A small delay in the inclusion [of your Note] in the *Comptes-rendus* will be much less detrimental to you than a premature publication of daring or insignificant results. Calm down a little."[28] To ease Villat's bitterness, Picard would later confide to him: "It is impossible to work at a leisurely pace today. As soon as one has written something, twenty mathematicians jump on it and one cannot let one's ideas grow more mature. [...] But you should not let this bother you too much: a house has a only one owner but a theorem can have many."[29] Pleased with the prospect of letting Brillouin handle Villat's "remarkable impetuosity" (Jacob 1979), Boussinesq left for the holiday and thought of something else.

---

[26] « Votre Note m'arrivait donc à point pour être présentée demain lundi, à la dernière séance de ce mois à laquelle je pourrai assister, devant partir vers le milieu de la semaine pour Les Vans (Ardèche). Mais je ne me sens pas en état de la présenter ; car le seul cas traité dans sa partie analytique conduit à la conclusion $P_x = 0$, justement contraire à votre conclusion principale et finale, relative, dîtes-vous, à un mouvement rotationnel permanent, tout-à-fait étranger à vos formules basées sur l'emploi d'un potentiel φ impliquant absence de tourbillons. Un note propre à motiver votre conclusion devrait ne pas contenir du tout de potentiel φ. » Boussinesq to Villat (11 September 1910). Villat Correspondence.

[27] « Votre Note n'apprendra donc rien à personne. » Boussinesq to Villat (11 September 1910). Villat Correspondance.

[28] « N'ayez, cher Monsieur, aucune crainte au sujet de la question de priorité : des travaux qui paraîtraient à l'étranger durant l'élaboration de votre thèse ne vous ôterait pas cette priorité à mes yeux, en tant que je pourrai être désigné pour juger cette thèse. Un peu de retard dans l'insertion aux Comptes-Rendus vous sera beaucoup moins nuisible que la publication prématurée de choses risquées ou insignifiantes. Calmez-vous donc. » Boussinesq to Villat (11 September, 1910). Villat Correspondance.

But Villat refused to calm down. On November 21, a short note of his was published in the *CRAS* (Villat 1910b). It was not Boussinesq, nor Brillouin who authorized its publication, but Émile Picard. Villat's work indeed seemed to be more interesting to mathematicians than to physicists and engineers. Early in 1911, he started to receive offers from various universities. A letter from the Lille mathematics professor Jean Clairin is instructive about the mathematical reception of Villat's work on fluid mechanics. Clairin had noticed Villat's articles in the *CRAS* and from them inferred that he was in possession of enough material for defending a doctoral thesis soon "and even an interesting thesis."[30] In the spring of 1911 Villat indeed rushed the defense of his thesis in order to get a position, not in Lille, but in Montpellier.[31] In a letter dated 16 February of that year, Picard agreed to write a short report on Villat's thesis since he was "in a hurry."[32] Picard's letter makes two points that reveal more about the interest mathematicians saw in Villat's work. First he noted that for a certain series

$$c_0 + c_1\zeta + c_2\zeta^2 + \cdots$$
$$+d_1\frac{1}{\zeta} + d_2\frac{1}{\zeta^2} + \cdots,$$

convergence on circumferences of radii $q$ and 1 was not obvious. Although they were banal for mathematicians, this type of reasoning was very foreign to applied mechanicians. Picard easily admitted that much, pointing out that this objection was irrelevant for applications, "but perhaps would there be a footnote to add."[33]

The second of Picard's remarks pointed out that Villat derived an integro-differential equation for a function $\Phi(\sigma)$ he had introduced that may be of interest for mathematicians (Villat 1911, 260). $\Phi(\sigma)$ was an angular function defined on the contour of an obstacle by means of which one could determine the form of

$$\Omega(t) = \frac{1}{\pi}\int_0^\pi \Phi(\sigma)\frac{1-t^2}{1-2t\cos\sigma+t^2}d\sigma.$$

---

[29] « On ne peut plus, hélas, travailler tranquillement aujourd'hui. Aussitôt que l'on a écrit quelque chose, vingt mathématiciens se jettent dessus, et on ne peut plus mûrir ses idées. […] Il ne faut pas cependant se préoccuper outre mesure ; une maison à un seul propriétaire, mais un théorème peut en avoir plusieurs. » Picard to Villat (16 February 1911).Villat Correspondence.

[30] « D'après les Notes de toi que j'ai vues récemment dans les Comptes rendus tu dois être en possession de résultats plus que suffisants pour une thèse — et même pour une thèse intéressante, dans ces conditions iln'est pas sans exemple qu'un candidat non docteur ait été préféré à un docteur. » Clairin to Villat (n.d. [1911]). Villat Correspondence.

[31] It is perhaps significant for Villat's decision to go to Montpellier that in the correspondence this university contrary to Lille appeared as one where the number of students was rather low. See Arnaud Danjoy to Villat (n.d. [1911]) and Clairin to Villat (n.d. [1911]). Villat Correspondence.

[32] « Puisque vous êtes pressé, je vais faire un rapport sommaire. » Picard to Villat (16 February, 1911). Villat Correspondence.

[33] « l'objection n'est pas bien importante pour les applications, puisqu'on peut se borner aux cas où elle n'aurait pas prise, mais peut-être cependant y aurait-il une remarque à ajouter. » Picard to Villat (16 February 1911). Villat Correspondence. A remark to that effect can indeed be found in (Villat 1911, 235).

This would lead Villat to the theory of integral equation, and Fredholm integrals especially, which was one of the liveliest research topics in analysis at the time. As it is well known, Villat's interest in such topic would develop over the years and be passed over to one of his most successful students, Jean Leray who would end up at the Collège de France after the Second World War.

## From Montpellier to Strasbourg: The Affirmation of a Mathematical Style

When Villat settled in Montpellier, his mathematical style was well established. To attack one of the most important technological problems of his days, he used the full extent of modern complex variables function theory and was flirting with the field of integro-differential equations, that is, with the forefront of mathematical research. His contemporaries were struck by his exceptional computation power which contrasted with the passion for set theory that has taken over some mathematicians of the time. "You are right," wrote his friend Léry:

> The best means to push forward mathematics is to eschew imprecision and generality […]. I am happy to see that you share this opinion, which is not the most common; since youth is now busy with sets, reasoning on hypotheses no one knows whether they might be realizable. [34]

Having read Villat's thesis, his friend Léry was dithyrambic: "I see that the supremacy of French mathematicians is not close to disappear." Written by—to say the least—a sympathetic critic that was himself a staunched nationalist, this assessment nevertheless seems significant to me.

Villat's style certainly was an expression of mathematical modernity. But contrary to the abstract and formalist modernism that was identified with the forefront of mathematics by the later Bourbakist generation, this one relied on precise and exact computation and on the audacious use of higher analysis for the solution of modern technological problems. This type of modernity was especially well represented in France where the usefulness of abstract analysis to mechanical applications was highly valued, as well as the impact technological concerns had on pure analysis itself. Such complementarities permeated contemporary discourses about the mathematician's social role. As the mechanician Gabriel Kœnigs wrote years later, Villat's "special esthetic" was the "happy answer of applications to theory."[35]

A widespread appreciation for Villat's mathematical is also evident from the rapid recognition his expository style—however forbidding—received. The dean of the Paris *faculté des sciences*, Paul Appell, was another specialist of higher analysis and rational mechanics who was enthused with Villat's work. He suggested his name to Jules Molk who was preparing the French adaptation of A. E. H. Love's article on hydrodynamics in Felix Klein's famous *Encyclopedia of the Mathematical Sciences* (Love, et al. 1912-1914). Molke asked him to review studies published between 1900 and 1912, adding: "there has not been a great number of them […], but there has been some and this you know better than I do."[36]

---

[34] « Tu as bien raison, le meilleur moyen de faire avancer les mathématiques est de ne pas rester dans le vague et le général […]. Je suis content de te voir cette opinion, qui n'est pas la plus répandue ; puisque la jeunesse fait des ensembles à tour de bras, en raisonnant sur des hypothèses dont on ne sait pas si elles sont réalisables. » Léry toVillat (n.d., 1911). « Tu as une puissance de calcul peu ordinaire. » Léry to Villat (24 January, 1911).Villat Correspondence.

[35] Comité secret du 14 mai 1923. Liste de présentation pour une place de Correspondant dans la Section de Mécanique. Rapport de M. KŒNIGS sur les travaux de M. Henri VILLAT, 8 p. Tapuscrit. Archives of the Academy of Sciences, Paris. Villat Biographical File.

Villat's mathematical style was further reaffirmed on the occasion of mathematical controversies that pitted against none other than Painlevé and Pierre Duhem. On 10 April, 1910, Painlevé addressed a long letter to Villat in which he objected to wake theory as a way to account for fluid resistance. Based on physical considerations, his argument was that without shocks and viscosity an infinite amount of energy was needed to establish the permanent regime wake theory took for granted. "It goes without saying, Painlevé wrote in concluding his letter, that the above remarks take nothing off the remarkable results you have obtained."[37] Only their practical value was put in question by Painlevé. This was of course irrelevant as far as Villat was concerned. Working in the best zero-viscosity limit imaginable, he was able to produced the best and most explicit approximation for resistance based on first principle. But that the limit may have been singular was well known to Villat.

Surprisingly enough, the second controversy in which Villat was involved occurred during the first months of World War I. In October 1914, Duhem used Umberto Cisotti's work to argue that even in the case where surfaces of discontinuity were introduced, the D'Alembert paradox still held (Duhem 1914a). In other words, supposing that the fluid was at rest at infinity, that the flow regime was permanent, and that an obstacle moved through it a at constant velocity, than no resistance to the motion could be due to the fluid. After publication in the *CRAS,* Levi-Civita, Picard and Villat wrote to Duhem to express their disagreement. Duhem's argument hinged on the consideration of integrals taken on a surface that was very far from the obstacle where, according to his hypothesis, the fluid was at rest. But as his correspondents pointed out in two-dimensional perfect fluid wake theory, lines of discontinuity extended to infinity. Duhem's integral therefore made no sense (Duhem 1914b, E. Picard 1914, Villat 1914).[38] In his letter to Villat, Picard reiterated the great esteem he had for Villat's work. He went on: "I have two sons and two sons-in-law in the battle; the last two are presently wounded. These are hours of anguish."[39]

This seems astonishing that Villat engaged in a minor controversy about integrals at infinity during this time of duress. In fact, little is known of Villat's activities during the war. His job file at Strasbourg states that he was mobilized from July 9, 1915 to September 24, 1918. Second-class soldier Villat was assigned to a ballistic center in Bourg d'Oisans where he computed range tables to fire against airships (Leray 1973). We also know that he fell ill, that he later had time to publish several papers, and that he received the prestigious Francœur prize from the Academy of Sciences. As for everyone in his generation, war is sure to have deeply affected him. His friend Léry died during the Marne battle on September 10, 1914. So had, two weeks earlier, on August 29, Jean Merlin, a mathematician by training working at the Lyon Observatory, whom Villat had befriended at the ÉNS.

By contrast, the direct effect war had on Villat's life no doubt were the tremendous opportunities the postwar period offered to him. In December 1918, the director of higher education for France, Alfred Coville, asked Villat if he were willing to go to Strasbourg on a temporary assignment to teach mechanics at the university.[40] Three months later, this assignment transformed into a chair with special advantages. After the Armistice, Strasbourg university was turned into a showcase of French academic excellence in German buildings. There, Villat had the rare advantage for a professor in the provinces of being able to teach advanced topics in fluid mechanics: modern hydrodynamics (1919), elliptic functions and their applications to fluid mechanics (1920), viscous fluids (1921), Laplace, Legendre, Fourier and Bessel functions and their applications to mathematical physics (1922), recent applications of analytical mechanics (1923), the theory of eddies (1924), wake theory (1925), and on the new hydrodynamics of low-viscosity fluids (1926). As several veterans went back to the university benches, he supervised some of the few doctoral theses to be defended outside of Paris: René Thiry, Maurice Roy, Jacques Berge, etc. (Leloup 2009).

---

[36] « il n'y en a pas eu un grand nombre depuis douze ans, mais enfin il y en a et vous le savez mieux que moi. » Molk to Villat (28 October, 1912).
[37] « Il va sans dire queles remarques précédentes n'enlèvent aucune valeur aux remarquables résultats mathématiques que vous avez obtenus. » Panlévé to Villat (10 April, 1913). Villat Correspondence.
[38] See also Duhem to Villat (2 November, 1914 ; 10 November, 1914) and Picard to Villat (13 November, 1914). Villat Correspondence.
[39] « J'ai deux fils et deux gendres à la bataille ; ces deux derniers sont actuellement blessés. Ce sont des heures d'angoisse. » Picard to Villat (13 November, 1914). Villat Correspondence.

As we can see, Villat's was heavily invested in the pedagogy of mathematical fluid mechanics. During the war, he wrote a long review of recent works on fluid mechanics (Villat 1917). In 1920, he published a small volume on fluid resistance which had become a crucial problem raised by wartime scientific mobilization due to its connection to both ballistics and aerodynamics (Villat 1920). Simultaneously, he was also asked to present his own work in Appell's *Traité de mécanique rationnelle*, 3rd ed., vol. 3 (Villat 1921a). As a result of this intense pedagogical activity, Villat was able to put his stamp on French mathematical fluid mechanics, a trend that was furthered in the late 1920s when he was appointed in Paris and published a series of influential textbooks.[41]

Before we come to that, let us note that Villat's situation in Strasbourg allowed him to exhibit his invaluable organization skills. In hastily organizing the controversial International Congress of Mathematicians in the symbolic site of Strasbourg in September 1920, Picard and the French mathematicians aimed at reorganizing the international community from which defeated nations were to be excluded (Lehto 1998, Lehto 2002). This would be the dawn of a "new order" (É. Picard 1920, 590). While Villat was left out of the political maneuverings, his material contribution did not go unnoticed. A he faced the difficulties of quickly organizing the congress in a country that had just come out of its greatest war, Picard wrote: "we are not opulent."[42] He had collected funds which he hoped would bring in 3000 francs. Villat took on the job and, in a few months, managed to collect 83,525 francs. The list of local donors published in the Proceedings of the Congress—of which Villat also was the editor—show the zeal with which Alsatian entrepreneurs contributed to it. As Kœnigs emphasized, Villat's "perseverant effort" and his "patriotic and scientific faith" contributed significantly to this financial success.[43]

Not only Villat's adminittrative skills, but his scientific authority was still admired by the men in place. Speaking of the mathematician Joseph Kampé de Fériet who had started to work on fluid mechanics because of the war, Brillouin wrote: "I hope that he also will get to these delicate questions; but until now he gives me the impression to limit himself to mathematical curiosities under the name and language of hydrodynamics instead of trying, like you have started to do, to forge the mathematical instrument needed to approach little by little the physical problem of hydrodynamics."[44] His scientific value coupled with his administrative efficiency thus made of Villat the right-hand man of several older "bosses" in the Parisian physico-mathematical community (Picard, Appell, Brillouin, Kœnigs) especially in academic publishing. From 1921, he would help Picard manage the *Annales scientifiques de l'École normale supérieure* as well as the *Bulletin des sciences mathématiques*. In 1922, he assumed the editorship of the *Nouvelles Annales de mathématiques*, a journal mostly geared at lycée professors of mathematics and, more importantly, the prestigious *Journal de mathématiques pures et appliquées* that had been founded by Joseph Liouville about eighty years earlier. In 1925, the Gauthier-Villars publishing house would ask him to head a series of high-level textbooks called the "Mémorial des sciences mathématiques," supposed to summarize the state of knowledge in a branch of mathematics. The first volume, by Appell, was titled *Sur une forme générale des équations de la dynamique*. By the mid-1920s, Villat had indeed become a central figure of the French mathematical establishment (Gispert et Leloup 2009).

---

[40] Coville to Villat (28 December, 1918 ; 20 March 1919). Villat Correspondence. On the University of Strasbourg, see (Crawford et Olff-Nathan 2005) and an article by Laurent Mazliak and Pavel Šišma in one of the upcoming collected volumes mentioned in note 5.

[41] One may also recall here the testimony of Paul Germain, a student of Villat's after WWII, who reported that he and his comrades were very impressed by the elegance of Villat's lectures and his great culture (Germain 2006, 107).

[42] Picard to Villat (28 February 1920). Villat Correspondence.

[43] The list of contributors is to be found in (Villat 1921b, xvi-xviii). Quotations from Kœnigs' speech are on p. xxvii.

# The Lure of Paris: The "Big Science" of Fluid Mechanics

The only missing pieces in order to become a central figure of French interwar mathematics now were for Villat to get a chair at the Sorbonne and a seat at the Academy of Sciences. Elected as correspondent of the Academy in 1923, he would get his seat in due time in 1932. But unfortunately for Villat, the situation in fluid mechanics was not looking good for him at the University of Paris. In this field, there were too many different people representing too many different interests and approaches. Suggesting that Villat accept a position as rector in Lille (which he did not take), Picard noted that this would not hinder his "coming back" to Paris. "I had hoped for a moment, he added, that the aeronautics business could be interesting for you but there are underpinnings in which I do not see clearly, despite [Jean] Villey's efforts."[45] As discussion went on in Paris, Picard reported a year later that a rather "violent" meeting of the Faculty Board had taken place apropos aeronautics: "very strong altercation between Kœnigs and Borel, Picard wrote. I wasn't there; it seems that [Lucien] Marchis also attacked Villey (who was not present). This matter has been spoiled, like so many things in our country where clans and chapels are proliferating."[46]

Since archival records of the Faculty Board at the Sorbonne are incomplete, the underpinnings of this dispute still have to be unearthed. At this time the scientific cluster built around aeronautical issues was taking a very innovative form at the University of Paris. Theoretical fluid mechanics teachings were complemented by courses on experimental hydrodynamics and aerodynamics and several conferences were held every year on current scientific or industrial topics. Marchis still occupied the aviation chair of the Sorbonne. At the head of the Saint-Cyr Aerotechnical Institute now was an engineer Albert Toussaint who had the status of associate professor [*maître de conferences*]. Having obtained, as war reparations, a series of publications from Prandtl's Institute in Göttingen and from the military aeronautics laboratory in Charlottenburg, Toussaint considerably revised experimental approaches to the problem of air resistance, but paid little attention to fundamental laws of motion (Toussaint 1921).

---

[44] « J'espère que lui aussi se mettra à l'œuvre sur ces questions délicates ; mais, jusqu'à présent il me fait l'impression de se borner à des curiosités mathématiques, sous le nom et le langage hydrodynamiques, au lieu de chercher, comme vous l'avez entrepris, à forger l'instrument mathématique nécessaire pour approcher peu à peu le problème physique de l'hydrodynamique ». Brillouin to Villat (9 September, 1919). Villat Correspondence.
[45] « si quelque chose de sérieux se présente pour vous à Paris, vous pourrez revenir. J'avais espéré un moment que l'affaire de l'aéronautique pourrait être intéressante pour vous, mais il y a là des dessous dans lesquels je ne vois pas bien clair, malgré les efforts de Villey qui voyait là quelque chose d'intéressant pour vous. » Picard to Villat (16 February, 1923). Villat Correspondence.
[46] « Il y a eu, il ya huit jours, une séance assez violente du Conseil de la Faculté [des Sciences de Paris], à propos de l'aéronautique : altercation très vive entre Kœnigs et Borel. Je n'y assistais pas ; il paraît que [Louis] Marchis a pris aussi Villey (qui n'était pas présent) à partie. Il y a làune affaire gâchée, image de bien des choses dans notre pays,où se multiplient les clans et les chapelles. » Picard to Villat (2 March 1924). Villat Correspondence.

After the war, the need for further scientific studies arose once more. In 1923, Painlevé was named to a chair in fluid mechanics that was established with governmental funds coming from the the aeronautic branch of the War Ministry [*sous-secrétariat d'état à l'aéronautique*]. But Painlevé barely taught there. He hired an engineer, Albert Métral, as associate professor and a Russian émigré, Dimitri Riabouchinski, to organize annual series of conferences at the Sorbonne, whose speakers included, over the years, Barillon, Bénard, Brillouin, Camichel, Hadamard, Maurice Roy, Yves Rocard, and in 1927, Villat himself.[47] In 1925 and 1926, Villat gave second-semester courses at the Sorbonne in Painlevé's chair.[48] Villat also taught at the *École nationale supérieure d'aéronautique*. In his course, he explains that to be concise, he excluded "tables of numbers and numerical results coming from detailed experiments."

After Painlevé's retirement, Villat was rather naturally named professor of fluid mechanics in January 1929. At this time, the situation had however evolved considerably. The political and technological context gave a tremendous boost to studies in fluid mechanics. In1928, the Air Ministry was created as a response to what was felt as a crisis in the air policy of the French state. One of the illustrious flying aces of the Great War, Charles de Kerland was enthusiastic about this new ministry. Underscoring that French motors were successfully exported, but not French gliders, he explained: "Our qualitative inferiority therefore is only partial and comes from the fact that we have not been able to conserve a perfect balance between all sides of the question" (Kerland 1929, 449-450). To take charge of scientific and technological affairs at the Air Ministry, a pioneer of aerostation was found who had also served as director of aeronautical research for the War Ministry in 1917–1918: Albert Caquot, whose initial assessment of the scientific situation was that: "Of all the branches of the scientific activity, one only had not taken, in France, the development required by the special need of aviation: […] Fluid Mechanics" (Caquot 1934, 4).

The program set up by Caquot may be counted among the biggest governmental effort for the development of a scientific discipline to date. A full-fledge Fluid Mechanics Institute [*Institut de mécanique des fluides*, herafter FMI] sponsored by the Air Ministry was established at the Sorbonne to be headed by Villat. The latter moreover found himself at the center of a network of special scientific institutes for fluid mechanics in Lille, Marseille, and Toulouse, as well as smaller training centers in Caen, Lyons, Nantes, Poitiers, and Strasbourg.[49] The neuralgic center of this wide network, the FMI of the University of Paris, kept a high level of funding throughout the Depression. It employed more than 20 people in theoretical and experimental researches, including Beghin, Bénard, and Foch who assured with Villat most of its teaching duties. Visiting professors now were sometimes coming from abroad, including Levi-Civita and Burgers. Experimental work was carried out in liaison with the Research Branch of the

---

[47] On fluid mechanics at the university of Paris and in France in the Interwar period, see (Fontanon et Franck 2006, 41-56, Anizan 2006, 129-134, 782-783, Mounier-Kuhn, L'enseignement supérieur, la recherche mathématique et la construction des calculateurs en France (1920-1970) 1998, 256).

[48] Two lessons in 1925 on "Leçons sur le mouvement d'un solide dans un fluide parfait, théorie des sillages." Four lessons in 1926 on "Le mouvement des solides dans un fluide peu visqueux. Les conditions limites quand la viscosité tend vers zéro. La nouvelle Hydrodynamique, théorie, applications et discussion."

[49] For a later historical survey of fluid mechanics institutes in France, see *Journées de mécanique des fluides (Marseille 1952),* "Publications scientifiques et techniques du ministère de l'Air," 296 (Paris, 1955). I want to thank Florian Schmalz for this reference.

Air Ministry: "We could not too much emphasize the support we received from that side" (Villat 1934, 426). The minister was apparently happy with the result, since on 10 October, 1933 he wrote to Villat: "the various courses of training set up by several universities […] will contribute more and more to direct researchers and scientists toward problems on which the progress of aeronautic industry relies" (Villat 1934, 441). As an historian has written, this was "big science" entering the mathematical sciences (Mounier-Kuhn 1998, 257-259).

Still, one is forced to acknowledge that the scientific outlook of those institutes, and of Villat's FMI in particular, seems at odds with contemporary research carried out in Germany or in the USA (Anderson 1997). The reigning paradigm for understanding fluid resistance remained Villat's wake theory, not Prandtl's boundary layer theory as envisioned by Toussaint. Engaged in a priority dispute with Theodore von Kármán about the eddie street appearing behind an obstacle in a fluid flow, Bénard had none of the technological acumen of the Hungarian engineer.[50] Only Riabouchinski, using his contacts among Russian émigrés' circles, was able to mobilized a large number of foreign scientists. Although opposition to the ministerial scheme to foster fluid mechanics research in France was rare, the sanguine physicist Henri Bouasse from the university of Toulouse expressed his skepticism at the mathematical approach suspecting the Minister of nepotism(Bouasse 1931, xx-xxiv). Although the archives of the institute were lost during World War II, the most notable legacy of the Paris FMI was the analog computers built by Joseph Pérès and Lucien Malavard before the war (Mounier-Kuhn 2010).

## Conflicted Ideologies About Applied Mathematics in Interwar France

The above sketch of how the study of fluid mechanics developed in interwar France leaves us with a paradoxical vision: a highly innovative complex associating State, university, the military, and private industry pursuing already antiquated scientific programs. Launched in the first decade of the century, these programs appear to have had little impact on contemporary and future technological development, despite initial intentions. A more nuanced appraisal of the development of research in fluid mechanics in France and the effects it may have had on industry during the interwar period will have to await further studies. For the moment, I would like to conclude this article by focusing on the conflicted ideologies of applied mathematics that arose from the experience of WWI. My claim is that this is useful to put Villat's rise to prominence in its proper context.

One may start by considering the official celebrations of the fiftieth anniversary of the French Mathematical Society, which took place on 24 May, 1924, at the Sorbonne. At the height of the modernist era in mathematics, the speeches published in the *Bulletin de la Société mathématique de France* on that occasion constitute what may seem like a surprising collection of statements. Welcoming his foreign colleagues (though still excluding German and Austrian mathematicians), Picard used a metaphor that would often be used by Bourbakists after WWII: "The mathematician's art is thus to create the molds in which physical theories try, at least for a time, to contain" the phenomena (SMF 1924, 29). After him, Lecornu discussed the applications of mathematics to mechanics, overconfidently pointing out that the problem of air resistance for aeronautic was "nearly solved [*à peu près résolue*]" (SMF 1924, 39). Others sketched the history of computing machines and the intimate link between mathematics and modern engineering. Borel painted a portrait of Henri Poincaré as embodying the encounter between mathematical speculation and practical life.

---

[50] On Bénard, see (Wesfreid 2006, Aubin 2008).

The Prime Minister Raymond Poincaré (Henri's cousin) soon joined in the concert of praise: "More than anyone else, Henri Poincaré has proved that the mathematical genius is not condemned to remain imprisoned in loneliness and abstraction, and that it is able to conceive wide outlooks over the world" (SMF 1924, 57).

Together these speeches composed a special portrait of mathematics where continuity with the past and interdependency between pure mathematics and neighboring disciplines. This picture, I want to claim, was heavily determined by the experiences of WWI. In a 1925 speech, the rector of the university of Paris for example declared that it had been "a true scientific epopee" (Lapie 1926, 45). From the start, WWI was widely perceived to have been a "scientific war" (see references in note 5), in which countless new applications of science were invented, developed, and used in battle. In ballistics and meteorology, for the development of new detection techniques such as sound ranging or the sonar, for the invention and development of new weapons such as the airplane and the tank, rather advanced mathematical tools had been mobilized by a great variety of actors. On 2 December, 1918, Painlevé thus exclaimed in front of the Academy of Science: "the most abstract and the most subtle mathematics partook in the solution of problems of detection [of enemy batteries] and in the computation of brand new range tables that increased the accuracy of fire by 25 percent" (Painlevé 1918, 809).

As far as mathematics was concerned, it would seem natural to think that the experiences of World War I as a scientific war may have produced regained interest for applied mathematics. In the opening speech he delivered at the Strasbourg ICM in 1920, Picard warned his colleagues against this perceived threat:

> Some say […] that applications of mathematics will be above all studies in the next years and that pure theory will be somewhat neglected by the generations. The times we are loving through have indeed become harder for the workers of intelligence, and the more optimists sometimes ask whether the civilization we are accustomed to will not be eclipsed. We therefore must not tire ourselves of repeating that in the last analysis the true source of all progress in the applied sciences lies in theoretical speculations (Villat 1921b, xxviii).

In his closing speech, Picard however also argued that since the world had completely changed between 1914 and 1920 and that the scientist now had to get out of his "ivory tower" (Villat 1921b, xxxii). The contradictory injunction—to resist utilitarianism by nurturing theoretical speculations while striving to be more involved in society—gave rise to several interesting professional trajectories in the postwar years. Injunctions addressed to mathematicians, and scientists in general, therefore were somewhat contradictory.

The contradiction is a striking feature of wartime speeches which extolled two extreme positions held by French Academicians with respect to German science. On the one hand, there were some who criticized the formal character of German scientific writing. In a pamphlet titled *L'histoire des sciences et les prétentions de la science allemande*, Picard was a chief proponent of the thesis according to which German scientists, impregnated with Kantian subjectivism, had produced another sort of scientific perversion detached from the reality of experience. "German science," he wrote in 1915, "has the tendency to posit *à priori* notions and concepts and to follow indefinitely the consequences, without worrying about their agreement with reality, and even while taking pleasure from distancing itself from common sense."[51] In an infamous book on *La Science allemande*, Duhem praised the realizations of German mathematicians but concurred with Picard in judging that from German science, excessive abstract "reasoning" at the expense of common sense "has banished Reason" (Duhem 1915, 22).

On the other hand, there was an opinion, perhaps more widespread, according to which it was utilitarianism and militarism that had perverted science in Germany. On 21 December, 1914, Appell used strong words to castigate the "barbarity" of the conception of science he saw at play among the France's enemies:

> The search for scientific truth […] is the noblest effort that a human being can ask from oneself. But, diverted from the constant ideal of right and humanity, going into the way of narrow specialization, disciplined in view of domination, [and ] mainly reduced go practical efficiency, the study of science quickly leads to a civilization characterized by its egotism, by its harshness, and by its materialism, to a kind of scientific barbarity [*barbarie savante*] like the one that has slowly taken over Germany today.[52]

After the initial shock, French Academicians realized that there might be a lesson to be learned from German science. Because of their alleged disinterest, great French savants like Pasteur or Berthelot had had the same success at turning their inventions into industrial products. On the opposite German utilitarianism excelled at it, Duhem further explained insisting on the dependency of all German invention on French scientific discoveries. In a time a war, this seemed like a crucial advantage that needed to change camp. Academic speeches were therefore inflected in this direction and the Germany's sense of organization began to be praised, even if the trope about German scientific barbarity remained.[53]

Today's historian may be prone to do away with the contradiction just pointed out simply by considering such speeches merely as the extravagant expression of exacerbated, irrational ideologies. But this would be a mistake. As we have seen, they were uttered by leaders of the postwar mathematical reconstruction who extensively recycled the themes they raised throughout the interwar period. In this context, Villat's highly mathematical approach to the problems of fluid mechanics seem more congenial. By showing that the full rigor of modern analytical theories may lead to results that were not unconnected to some of the practical problems arising from modern technologies, Villat seemed to be able to navigate between the Scylla of senseless abstraction and the Charybdis of utilitarianism.

This view however may not have been wholly convincing to the practitioners. A physics textbook about *Résistance des fluides* was for example prefaced with a rabid denunciation of the abuse of deductive reasoning in the physical sciences (Bouasse, Résistance des fluides: vol des avions et des oiseaux, hélices et moulins à vent, manœuvre des navires 1928). The usefulness of mathematics for the training of physicists and engineers was a strongly debated issue before and after WWI, but during the war this debate had taken on a highly ideological charge.[54] In 1923, Bouasse had chosen not to open his treatise on hydrodynamics by general equations of motion for "two excellent reasons: to be understandable and not to bore. General equations are mainly used by mathematicians to show that they are very skilful, but also that they have little concern for the phenomena" (Bouasse 1923). Two years earlier, he already prefaced another textbook by an argument "On the Uselessness of Mathematics for the Training of the Mind [*esprit*]." While Bouasse rarely straightforwardly attacked Villat's approach, the reviews of his books provided an occasion for L. Potin to expressed more direct criticisms: "if manufacturers had counted on conformal representation to improve their materials, he wrote, airplanes would not have achieved their present perfection."[55] In a later review, Potin went as far as claiming that discontinuous stream theory had been "abandoned [since] the introduction of surfaces of discontinuity was deemed arbitrary."[56]

---

[51] *CRAS* 161 (1915), 411.
[52] *CRAS* 159 (1914), 822.
[53] See for example the geologist Edmond Perrier's speech in *CRAS* 163 (1916), 21-22, Lecornu's in ibid. 166 (1918), 407, or Painlevé's in ibid. 167 (1918), 800.
[54] On the impact on this debate on teachings at the École poytechnique, see an article by Jean-Luc Chabert and Christian Gilain in one of the upcoming collected volumes mentioned in note 5.

The frequency of such judgments may help to place public discourses about mathematics, rigor, and applications in their proper context. The middle position may have been best expressed by Robert d'Adhémar, an engineer and self-taught mathematician who had defended a thesis in front of Picard in 1904. Reacting to Bouasse's provocative preface from 1921, he wrote: "there is a latent conflict between physicists and mathematicians." It is a "perversion," both authors agreed, to identify rational mechanics with analysis; instead, mechanics should be taught as the first chapter of physics courses by insisting on practice and numerical results. But unlike Bouasse Adhémar refused to conclude that a strong mathematical culture was useless. He defended the rigor of mathematicians' reasoning and insisted on the need to adapt mathematics in order to make them useful. "War," Adhémar concluded, "showered us with realism": analysis, mechanics and physics should be taught in a more integrated manner (Adhémar 1921, 275). In this context, the special appeal of Villat's highly mathematical approach to fluid mechanics seemed more obvious. Indeed, Adhémar also reviewed Villat's first book on fluid resistance (Villat 1920) in the following terms:

> I must warn the reader that only the mathematician can read M. Villat. […] When a question is difficult and requires knowledge and skill, the brainless have the resources […] to say: "This is theory." Let us congratulate those who like M. Villat make beautiful theories which are the necessary framework of science in the making.[57]

## Conclusion

At the 1924 meeting celebrating the SMF's anniversary at the Sorbonne, the Belgian mathematician Charles de la Vallée-Poussin was, as the acting president of the International Mathematical Union, the only foreigner to address the assembly. Echoing the most extreme views uttered during the last war by some of his audience, he emphasized that science knew no border: "there is neither German nor French science, […] there is just Science" (SMF 1924, 34). This was true, he acknowledged. But only in some sense, he went on:

> because mathematical science is an art and because an art leaves a wide room to individual sensitivity and esthetics […], taste that is the product of education, of the milieu of some deep-rooted sub-conscience coming from the genius of the race (SMF 1924, 34).

De la Vallée-Poussin underscored that for the French taste rigorous abstraction "smell of the laboratory" (SMF 1924, 35). Mathematics, he went on, owed to France more than the important contributions made by its mathematicians, but something coming from the spirit of this nation: "mathematicians […] like Frenchmen […] possessed the art of filtering murky notions and the gift of clarifying obscure ideas, most particularly those of the Germans" (SMF 1924, 35-36). The racist tone of this speech clearly shows the persistence of ideas expressed in the special circumstances of the war. This was a potent public discourse that was wholly coherent with extreme wartime views according to which senseless abstraction was as dangerous morally speaking as utilitarianism. From the violent condemnation of what was perceived as "German" perversions and barbarity, a new consensus about the nature and goals of mathematics had emerged in France, whereby a strict adhesion to the full rigor of higher analysis could be wed with a concern for applications that did not translate into immediate applicability. While this may have seemed absurd to the practitioners, this was a posture that swayed many in the mathematical establishment.

---

[55] *RGSPA* 39 (1928), 312.
[56] *RGSPA* 39 (1928), 343.
[57] *RGSPA* 33 (1922), 152.

The thesis I have defend here may be conceived as a French counterpart to Paul Forman's famous argument. By claiming that German mathematicians and physicists rejected causality and embraced a certain form of modernity as a result of their adapting to a hostile intellectual environment (Forman 1971), the historian drew a lot of criticism. My argument however does not concern strong commitment to a hypothesis about the nature of fluid flows, but rather the privileging of a specific approach to deal with some of the scientific problems underlying aeronautics. As such, it shows that close attention to public discourses can indeed help understand the contours taken by a national scientific community in the years that followed WWI.

I have argued that the paradoxes Villat's career dissolve in this context. Propelled by a pressing technological problem of modernity, air flight and the air resistance problem, wake theory revealed the profound beauty of modernist mathematical methods. In his report on the newly established Fluid Mechanics Institute of the University of Paris, Villat claimed that the effort of the Air Ministry were warranted solely on the ground that due to the courses and laboratories organized by the Ministry, one of his students, Jean Leray, was led to the abstract mathematical work for which he has become famous.[58] This perhaps is the ultimate paradoxical twist in the story: that for Villat the justification of the unprecedented effort of a French government to promote science in view of developing the aeronautical industry was an abstract mathematical theory that was held in great consideration by the Bourbakist generation!

## Abbreviations
BSMF : *Bulletin de la Société mathématique de France*
CRAS : *Comptes rendus de l'Académie des sciences*
ÉNS: École normale supérieure, Paris.
FMI : Fluid Mechanics Institute, University of Paris.
RGSPA : *Revue générale des sciences pures et appliquées.*
SMF : the French Mathematical Society [*Société mathématique de France*].
WWI : World War I.

---

[58] Known axxx, s Leray's work directly inspired David Ruelle's work at the root of chaos theory as late as 1970. Leray's justifying the ffort of the Ministry is to be foundin (Villat 1934, 429). This is also confirmed by Leray to Villat (n.d. [1934]). Villat Correspondence.

# Bibliography

Adhémar, Robert d'. "Physiciens et mathématiciens." *Revue générale des sciences pures et appliquées*, 1921: 273-275.

Anderson, John David. *A History of Aerodynamics and Its Impact on Flying Machines.* Cambridge: Cambridge University Press, 1997.

Anizan, Anne-Laure. *Paul Painlevé (1863-1933): un scientifique en politique.* Doctoral thesis, Paris: Institut d'études politiques, 2006.

Anonymous. "Progress of Mechanical Flight." *Flight*, 1909: 12.

Armengaud, Jules. "Le problème de l'aviation et sa solution par l'aéroplane." *Bulletin de la Société d'encouragement pour l'industrie nationale*, juillet 1907: 847-873.

Aubin, David. "'The Memory of Life Itself': Bénard's Cells and the Cinematograph of Self-Organization." *Studies in the History and Philosophy of Science*, 2008: 359-369.

Aubin, David, and Patrice Bret. *Le Sabre et l'éprouvette: l'invention d'une science de guerre, 1914-1939.* Paris: Agnès Viénot, 2003.

Beaulieu, Liliane. "Regards sur les mathématiques françaises entre les deux guerres." *Revue d'histoire des sciences*, 2009: 9-37.

Biot, Maurice Anthony. "L'hydrodynamique moderne et ses applications." *Revue des questions scientifiques*, September 1930: 235-257.

Bouasse, Henri. *Jets, tubes et canaux.* Paris: Delagrave, 1923.

—. *Résistance des fluides: vol des avions et des oiseaux, hélices et moulins à vent, manœuvre des navires.* Paris: Delagrave, 1928.

—. *Résistance des fluides: vol des avions et des oiseaux, hélices et moulins à vent, manoeuvre des navires.* Paris: Delagrave, 1928.

—. *Tourbillons, forces acoustiques, circulations diverses.* Paris: Delagrave, 1931.

Brillouin, Marcel. "Les surfaces de glissement d'Helmholtz et la résistance des fluides." *Annales de chimie et de physique* 23 (1911): 145-230.

Brillouin, Marcel. "Mouvement discontinu de Helmholtz. Obstacles courbes." *CRAS* 151 (1910): 931-933.

Caquot, Albert. *Notice complémentaire sur les travaux historiques.* Paris : Guthier-Villars, 1934.

Carlier, Claude. "Ferdinand Berber et l'aviation." *Guerres mondiales et conflits contemporains*, 2003: 7-23.

Chadeau, Emmanuel. "État, industrie, nation: la formation des technologies aéronautiques en France (1900-1950)." *Histoire, économie et société*, 1985: 275-300.

Champonnois, Sylvain. "Les Wright et l'armée française: les débuts de l'aviation militaire (1900-1909)." *Revue historique des armées*, 2009: 108-121.

Crawford, Elisabeth, and Josiane Olff-Nathan. *La Science sous influence: l'université de Strasbourg, enjeu des conflits franco-allemands 1872-1945.* Strasbourg: La Nuée bleue, 2005.

Darrigol, Olivier. *Worlds of Flow: A History of Hydrodynamics from the Bernoullis to Prandtl.* Oxford: Oxford Univ. press., 2005.

Didi-Hubermann, Georges, and Laurent Mannoni. *Mouvements de l'air: Étienne-Jules Marey, photographe des fluides.* Paris: Gallimard, 2004.

Duhem, Pierre. "Remarque sur le paradoxe hydrodynamique de d'Alembert." *Comptes-rendus de l'Académie des sciences*, 1914b: 638-640.

—. "Sur le paradoxe de Dalembert." *Comptes-rendus de l'Académie des sciences*, 1914a: 592-595.

—. *La Science allemande.* Paris: A. Hermann, 1915.

Eiffel, Gustave. *Recherches expérimentales sur la résistance de l'air exécutées à la tour Eiffel.* Paris: L. Maretheux, 1907.

Fontanon, Claudine. "La naissance de l'aérodynamique expérimentale etses applications à l'aviation.Une nouvelle configuration socio-technique." In *Histoire de la mécanique appliquée: Enseignement, recherche et pratiques mécaniciennes en France après1880*, by Claudine Fontanon, 57-88. Paris: SFHST & ENS édidtions, 1998.

Fontanon, Claudine, and Robert Franck. *Paul Painlevé: Un savant en politique.* Paris: Presses universitaires de Rennes, 2006.

Forman, Paul. "Weimar Culture, Causality, and Quantum Theory, 1918-1927: Adaptation by German Physicists and Mathematicians to a Hostile Intellectual Environment." *Historical Studies in the Physical Sciences*, 1971: 1-115.

Germain, Paul. *Mémoires d'un scientifique chrétien* . Paris: L'Harmatthan, 2006.

Gispert, Hélène, and Juliette Leloup. "Des patrons des mathématiques en France dans l'entre-deux-guerres." *Revue d'histoire des sciences*, 2009: 39-118.

Goldstein, Catnerine. "La théorie des nombres en France dans l'entre-deux-guerres: de quelques effets de la première guerre mondiale." *Revue d'histoire des sciences*, 2009: 143-175.

Hadamard, Jacques. "Lery (Georges), né à Limours le 28 avril 1880, tué à l'ennemi le 10 septembre 1914. Promotion de 1899." *Association amicale de secours des anciens élèves de l'Ecole normale supérieure (Paris), Réunion générale annuelle*, 1916: 112-116.

Hashimoto, Takehiko. "Theory, Experiment, and Design Practice: The Formation of Aeronautical Research, 1909-1930." *Ph. D. dissertation.* Baltimore: Johns Hopkins University, 1990.

Henri.

Jacob, Caius. "Henri Villat, membre d'honneur de l'Académie Roumaine, sa vie et son œuvre." *Noesis : travaux du comité roumain d'histoire et de philosophie des sciences*, 1979: ???

Kerland, Charles Le Coq de. "Le Ministère de l'Air." *Revue politique et parlementaire*, 1929: 447-475.

Lapie, Paul. "Discours du recteur." *Annales de l'université de Paris*, 1926: 37-47.

Lecornu, Léon. "Revue de mécanique appliquée." *Revue générale des sciences pures et appliquées*, 1909: 30-42.

Lee, James Lawrence, Claudine Fontanon, and Daniel H. Fruman. *Eiffel 1903 Drop Test Machine and 1912 Wind Tunnel.* Paris: Historical Mechanical Engineering Landmarks, 2005.

Lehto, Olli. *Mathematics without Borders: A History of the International Mathematical Union.* New York: Springer, 1998.

Lehto, Olli. "The Formation of the International Mathematical Union." In *Mathematics Unbound: The Evolution of an International Mathematical Research Community, 1800-1945*, by Karen Hunger Parshall and Adrian Rice, 381-395. Providence: American Mathematical Society, 2002.

Leloup, Juliette. "L'entre-deux-guerres mathématique à travers les thèses soutenues en France." *doctoral thesis.* Paris: Université Pierre et Marie Curie, 2009.

Leray, Jean. "Villat (Henri), né à paris le 24 décembre 1879, décédé à Paris le 19 mars 1972.—Promotion de 1899." *Association amicale de secours des anciens élèves de l'Ecole normale supérieure (Paris), Réunion générale annuelle*, 1973: 37-39.

Levi-Cività, Tullio. "Sci e leggi di resistenza." *Rendiconti del Circolo matematico di Palermo* 23 (1907): 1-37.

Love, A. E. H., Paul Appell, Henri Begin, and Henri Villat. *Développements concernant l'hydrodynamique.* Vols. IV-5, in *Encyclopédie des sciences mathématiques pures et appliquées*, by Jules Molk, Paul Appell, Felix Klein and C. H. Müller. Paris: Gauthier-Villars, 1912-1914.

Maurain, Charles. "Revue d'aérotechnique expérimentale." *Revue générale des sciences pure et appliquées*, 1914: 679-686.

—. "Les études à l'Institut aérotechnique de Saint-Cyr." *Revue générale des sciences pures et appliquées*, 1913: 218-229.

Mazliak, Laurent, and Rossana Tazzioli. *Mathematicians at Wat: Volterra and His French Colleagues in World War I.* Berlin: Springer, 2010.

Mounier-Kuhn, Pierre-Éric. "L'enseignement supérieur, la recherche mathématique et la construction des calculateurs en France (1920-1970)." In *Des ingénieurs pour la Lorraine: enseignements industriels et formations technico-scientifiques supérieures, XIXe-XXe siècles*, by F. Birck and André Grelon, 251-286. Metz: Éd. Sepernoise, 1998.

—. *L'informatique en France de la Seconde Guerre mondiale au Plan Calcul: l'émergence d'une science.* Paris: Presses de l'université Paris-Sorbonne, 2010.

Nastasi, Pietro, and Rossana Tazzioli. "Toward a Scientific and Personal Biographyof Tullio Levi-Civita (1873-1941)." *Historia Mathematica*, 2005: 203-236.

Painlevé, Paul. "Séance publique annuelle." *Comptes rendus de l'Académie des sciences*, 1918: 797-810.

Painlevé, Paul. "L'aéroplane." In *Compte rendu de la 38e session de l'Association française pour l'avancement des sciences à Lille*, i-xii. Paris: AFAS, 1910.

Picard, Emile. "A propos du paradoxe hydrodynamique de d'Alembert." *Comptes-rendus de l'Académie des sciences*, 1914: 638.

Picard, Émile. "Le Congrès international de Mathématiques de Strasbourg." *Comptes-rendus de l'Académie des sciences*, 1920: 589-591.

Prochasson, Christophe, and Anne Rassmussen. *Au nom de la patrie: les intellectuels et la Première Guerre mondiale (1910-1919).* Paris: La Découverte, 1996.

Renard, Paul. "La Résistance de l'air et les récentes expériences de M. G. Eiffel." *Revue générale des sciences pures et appliquées*, 1909: 78-84.

Ricca, Renzo L. "The Contributions of Da Rios and Levi-Civita to Asymptotic Potential Theory and Vortex Filament Theory." *Fluid Dynamic Research* 18 (1996): 245-268.

Robert. "Les recherches de physique dans leur rapport avec l'aérotechnique." In *Le Livre du cinquantenaire de la Société française de physique*, 73-97. Paris: Éditions de la Revue d'optique théorique et expérimentale, 1925.

Roussel, Yves. "Histoire d'une politique des inventions, 1887-1918." *Cahiers pour l'histoire du CNRS*, 1989: 19-57.

Roy, Maurice. «Notice nécrologique sur Henri Villat, membre de la section de mécanique.» *Comptes-rendus de l'Académie des sciences*, 1972: 127-132.

SMF. "Vie de la société." *Bulletin de la société mathématique de France*, 1924: 1-67.

Soreau, Rodolphe. "Review of Notes sur l'aérodynamique de l'aéroplane, par Emmanuel Vaillier." *Mémoires de la Société des ingénieurs civis*, 1905: 820-821.

Toussaint, Albert. "Revue d'aérotechnique expérimentale." *Revue générale des sciences pures et appliquées*, 1921: 203-213.

Villat, Henri. "Institut de Mécanique: son activité en 1933. I - Mecanique des fluides." *Annales de l'Université de Paris*, 1934: 425-441.

—. "Sur la résistance des fluides." *Annales scientifiques de l'École normale supérieure*, 1911: 203-311.

—. *Aperçus théoriques sur la résistance des fluides.* Paris: Gauthier-Villars, 1920.

—. "Quelques récents progrès des théories hydrodynamiques." *Bulletin des sciences mathématiques*, 1917: 43-60.

—. *Comptes rendus du Congrès international des mathématiciens (Strasbourg, 22-30 septembre 1920).* Toulouse: Édouard Privat, 1921b.

—. *Comptes-rendus de l'Académie des sciences*, 1914: 800-803.

Villat, Henri. *Mouvement plan discontinu d'un liquide.* Vol. 3, in *Traité de mécanique rationnelle*, by Paul Appell, 523-561. Paris: Gauthier-Villars, 1921a.

—. "Sur le théorème de Lagrange en hydrodynamique." *Nouvelles annales de mathématiques*, 1910a: 282-284.

Villat, Henri. "Sur la résistance des fluides limites par une paroi fixe indéfinie." *CRAS* 151 (1910b): 933-935.

Villat, Henri. "Sur les mouvements d un fluide autour d'un obstacle de forme donnée." *CRAS* 151 (1910c): 1034.

Wesfreid, J. Eduardo. "Scientific Biography of Henri Bénard (1874–1939)." In *Dynamics of spatio-temporal cellular structures: Henri Bénard centenary review*, by I. Mutabazi, J. E. Wesfreid and É. Guyon, 9-40. New York: Springer, 2006.

# Figures :

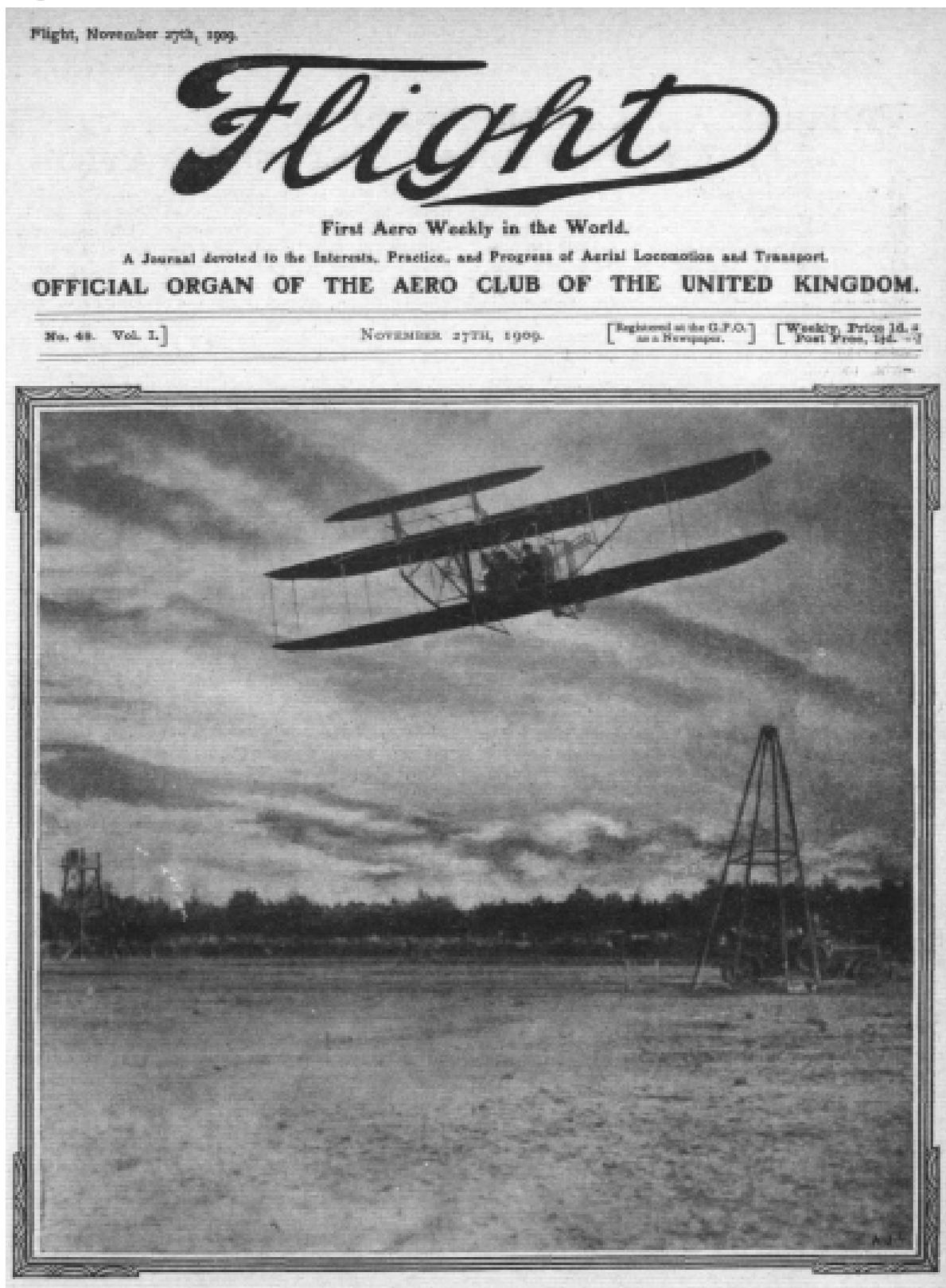

**Figure 1 :** Front page of *Light* (27 November 109) showing the historic flight of Paul Painlevé with Wilbur Wright at Auvours, France, on 11 October, 1908.

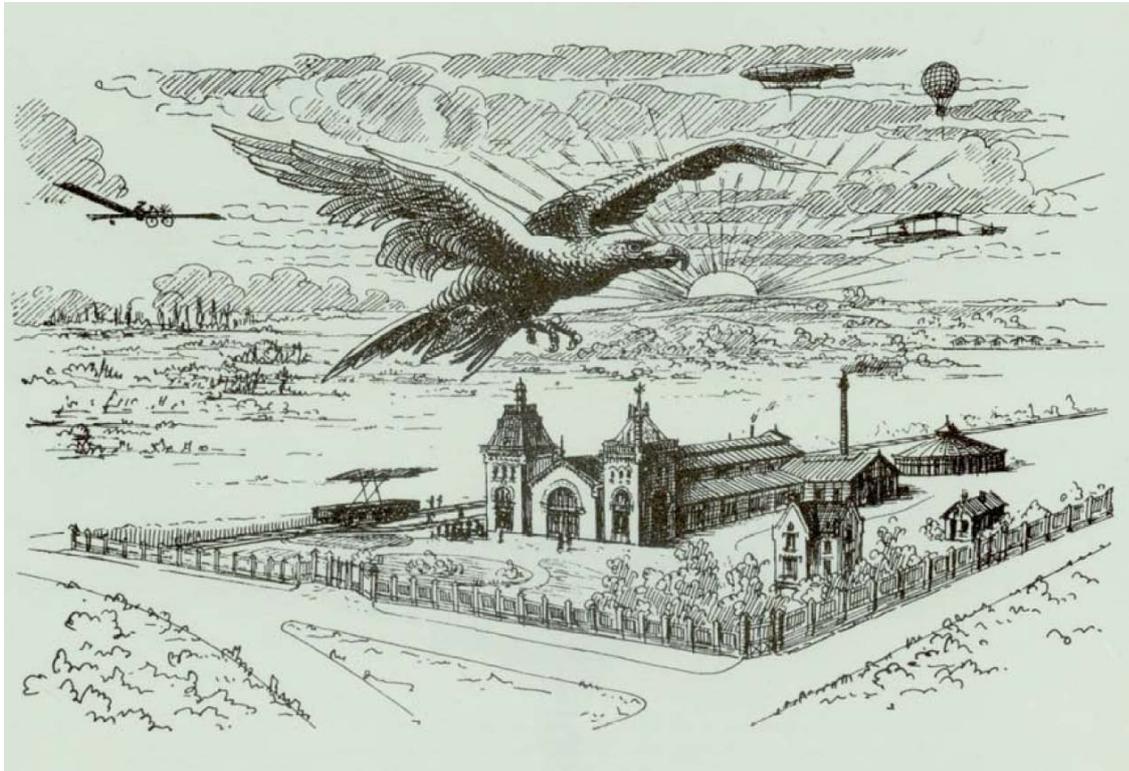

**Figure 2 :** Institut aérotechnique de l'université de Paris in Saint-Cyr. *Bulletin de l'Institut aérotechnique de l'université de Paris 1 (1911).*

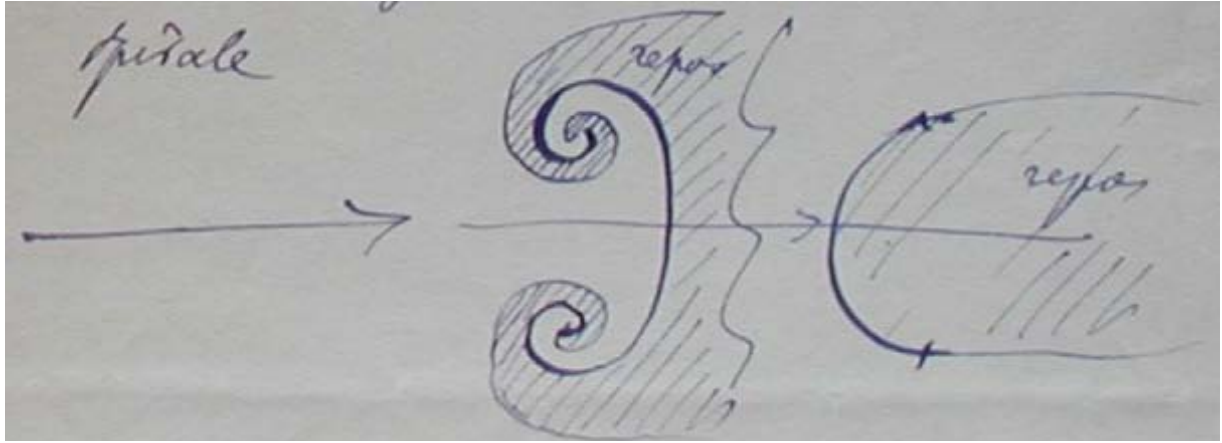

**Figure 3 :** Examples of surfaces treated by Marcel Brillouin for the problem of perfect fluid flows around an obstacle in his letter to Villat, 27 July 1910. © Reproduced with the kind permission of the Academy of Science in Paris.

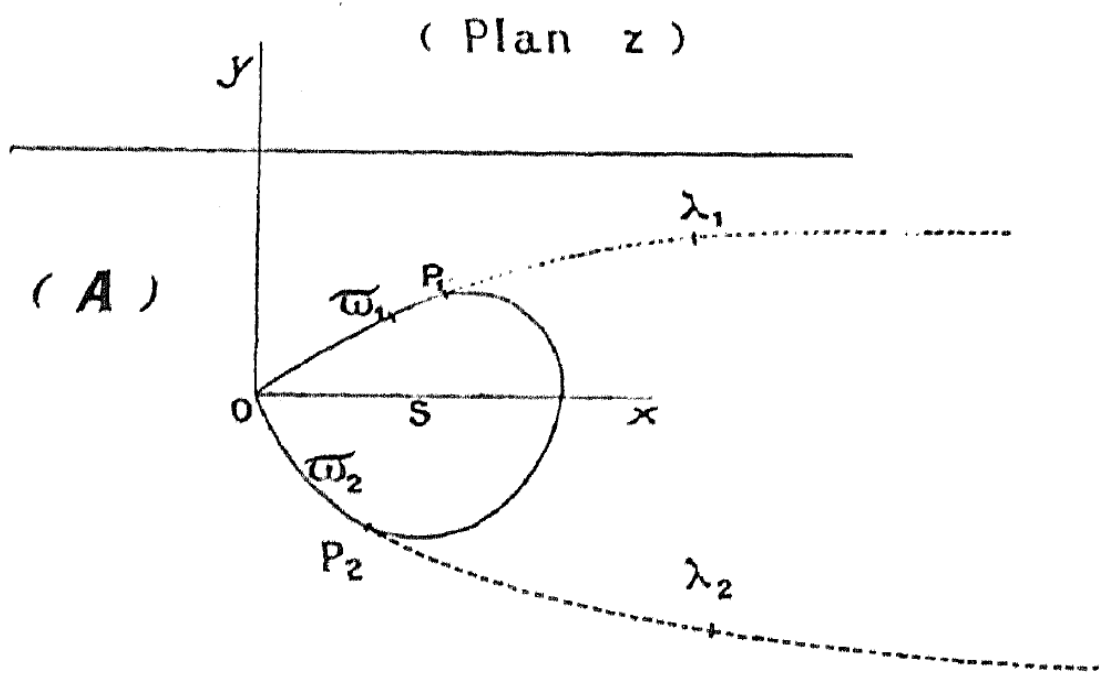

**Figure 4 :** Representation of the obstacle in the *z*-plane corresponding to the physical motion of the fluid around the obstacle (**Villat 1911, 209**).

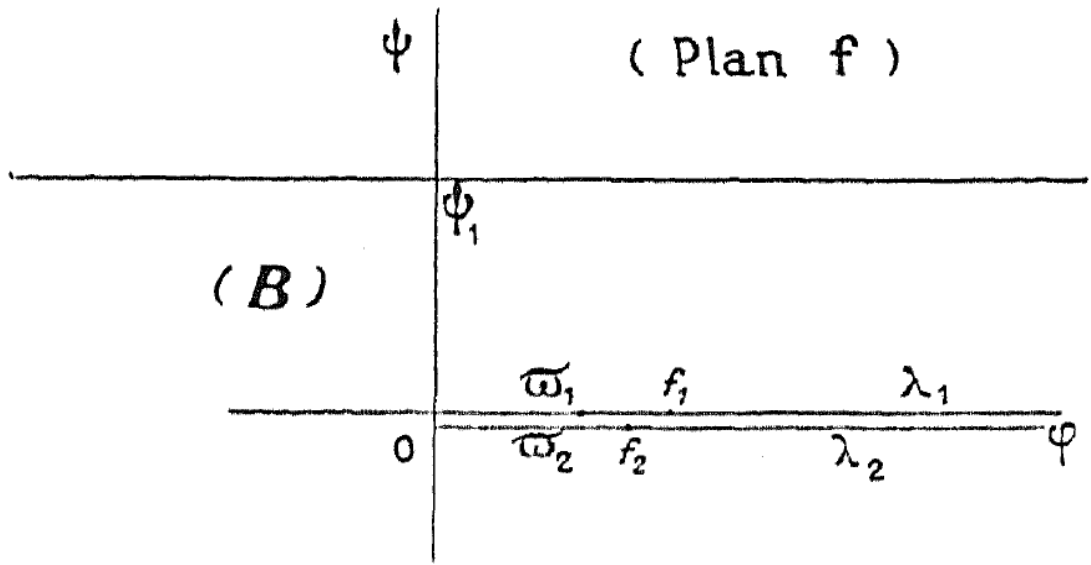

**Figure 5 :** Representation of the fluid motion in the *f*-plane (**Villat 1911, 211**).

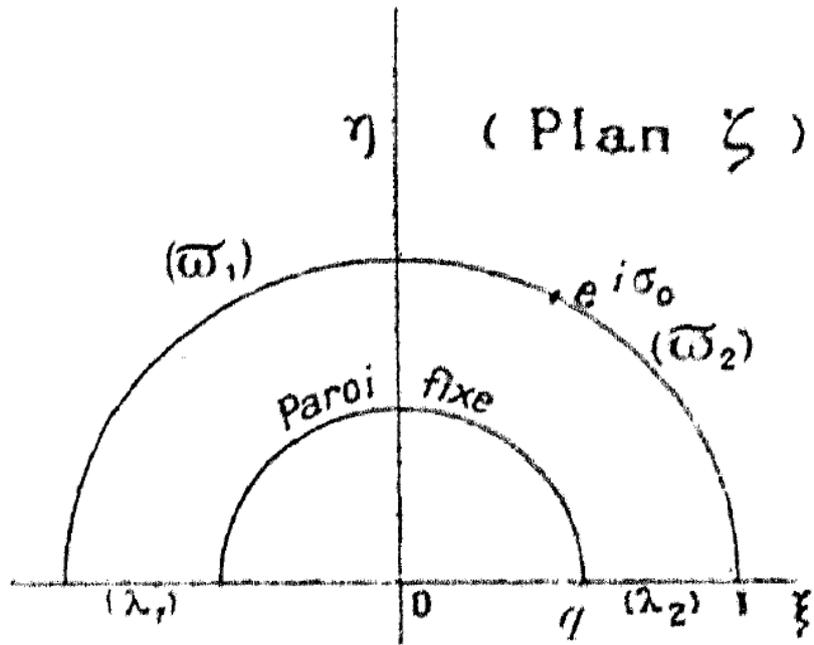

**Figure 6 :** Representation of the fluid motion in the ζ-plane (**Villat 1911, 216**).

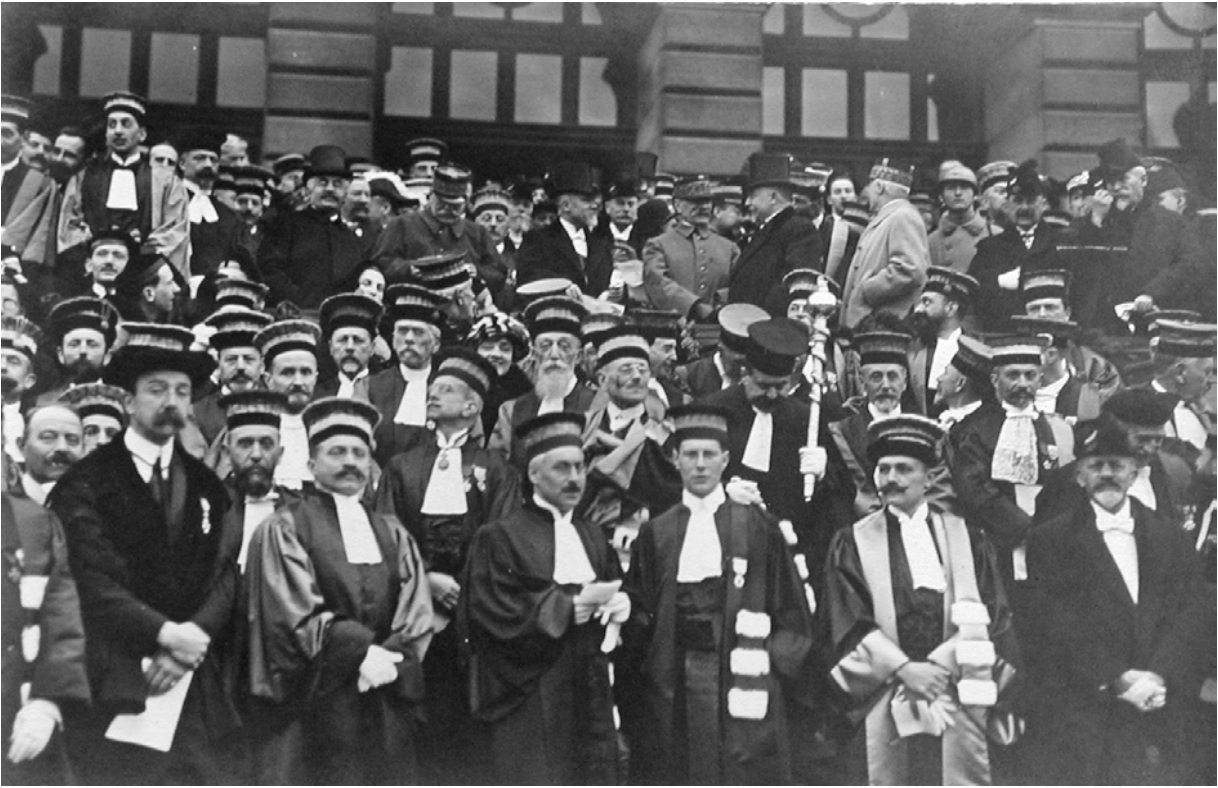

**Figure 7:** Officials Officers, and Professors in Front of the Palais universitaire in Strasbourg, 1919. Courtesy from the Archives del'Académie des Sciences in Paris.